\documentclass[a4paper, 11pt]{article}
\pdfoutput=1 

\usepackage[T1]{fontenc} 
\usepackage[utf8]{inputenc} 

\usepackage{amsfonts,amsmath,amssymb,mathtools,dsfont} 
\usepackage{microtype} 
\usepackage[usenames, dvipsnames]{xcolor} 

\usepackage{fancyhdr}

\usepackage[labelsep = period, labelfont = bf, justification = centering]{caption}
\usepackage{float,graphicx,subcaption} 

\usepackage{booktabs,makecell} 
\usepackage{enumitem,mdwlist} 
\setlist[itemize]{topsep=0ex,itemsep=0ex,parsep=0.4ex}
\setlist[enumerate]{topsep=0ex,itemsep=0ex,parsep=0.4ex}

\usepackage{etoolbox}
\usepackage[super]{nth} 
\allowdisplaybreaks

\usepackage[left = 2cm, right = 2cm, top = 2.5cm, bottom = 2.5cm, headsep = 12pt, headheight = 15pt]{geometry}
\usepackage{parskip}

\usepackage[noBBpl]{mathpazo} 
\DeclareFontFamily{U}{matha}{\hyphenchar\font45}
\DeclareFontShape{U}{matha}{m}{n}{
	<5> <6> <7> <8> <9> <10> gen * matha
	<10.95> matha10 <12> <14.4> <17.28> <20.74> <24.88> matha12
}{}
\DeclareSymbolFont{matha}{U}{matha}{m}{n}

\DeclareMathSymbol{\specialuparrow}{\mathrel}{matha}{"D2}
\DeclareMathSymbol{\specialrightarrow}{\mathrel}{matha}{"D1}

\DeclareFontFamily{U} {cmmi}{}

\DeclareFontShape{U}{cmmi}{m}{n}{
	<-6> cmmi5
	<6-7> cmmi6
	<7-8> cmmi7
	<8-9> cmmi8
	<9-10> cmmi9
	<10-12> cmmi10
	<12-> cmmi12}{}

\DeclareSymbolFont{Xcmmi} {U} {cmmi}{m}{n}

\DeclareMathSymbol{\mu}{\mathord}{Xcmmi}{'026}
\DeclareMathSymbol{\rho}{\mathord}{Xcmmi}{'032}
\DeclareMathSymbol{\varphi}{\mathord}{Xcmmi}{'047}

\DeclareFontFamily{U} {cmr}{}

\DeclareFontShape{U}{cmr}{m}{n}{
	<-6> cmr5
	<6-7> cmr6
	<7-8> cmr7
	<8-9> cmr8
	<9-10> cmr9
	<10-12> cmr10
	<12-> cmr12}{}

\DeclareSymbolFont{Xcmr} {U} {cmr}{m}{n}
\DeclareMathSymbol{\Delta}{\mathord}{Xcmr}{'001}
\DeclareMathSymbol{\Upsilon}{\mathord}{Xcmr}{'007}
\DeclareMathSymbol{\Omega}{\mathord}{Xcmr}{'012}

\usepackage[hyphens]{url} 
\usepackage[linktoc = all, hidelinks, colorlinks, unicode=true, pagebackref=true]{hyperref} 
\usepackage[capitalise, compress, nameinlink, noabbrev]{cleveref} 

\hypersetup{linkcolor={blue!70!black}, citecolor={green!70!black}, urlcolor={blue!70!black}}

\crefname{section}{\S}{\S\S} 
\Crefname{section}{Section}{Sections} 
\crefname{subsection}{\S}{\S\S} 
\Crefname{subsection}{Subsection}{Subsections} 
\crefname{page}{page}{pages}
\Crefname{page}{Page}{Pages}

\usepackage{pgf,tikz}
\usepackage{mathrsfs}
\usetikzlibrary{decorations.pathreplacing}
\usetikzlibrary[patterns]
\usetikzlibrary{arrows.meta}
\usepackage{tkz-euclide}

\tikzset{font={\fontsize{10pt}{12}\selectfont}}
\tkzSetUpPoint[fill=black, size = 3pt]
\tkzSetUpLine[color=black, line width=0.6pt]

\usepackage{amsthm,thmtools,thm-restate}

\declaretheoremstyle[
spaceabove = .5\baselineskip\@plus.2\baselineskip\@minus.2\baselineskip, 
spacebelow = .2\baselineskip\@plus.2\baselineskip\@minus.2\baselineskip,
headfont = \normalfont\bfseries,
notefont = \mdseries, 
notebraces = {(}{)},
bodyfont = \normalfont\itshape,
postheadspace = .5em,
headpunct = .
]{bolditalic}

\declaretheoremstyle[
spaceabove = .5\baselineskip\@plus.2\baselineskip\@minus.2\baselineskip, 
spacebelow = .2\baselineskip\@plus.2\baselineskip\@minus.2\baselineskip,
headfont = \normalfont\bfseries,
notefont = \mdseries, 
notebraces = {(}{)},
bodyfont = \normalfont,
postheadspace = .5em,
headpunct = .
]{boldnormal}

\declaretheoremstyle[
spaceabove = .2\baselineskip\@plus.2\baselineskip\@minus.2\baselineskip, 
spacebelow = .5\baselineskip\@plus.2\baselineskip\@minus.2\baselineskip,
headfont = \normalfont\itshape,
notefont = \mdseries, 
notebraces = {}{},
bodyfont = \normalfont,
postheadspace = .5em,
headpunct = .,
qed = \qedsymbol
]{proofstyle}

\declaretheoremstyle[
spaceabove = .5\baselineskip\@plus.2\baselineskip\@minus.2\baselineskip, 
spacebelow = .2\baselineskip\@plus.2\baselineskip\@minus.2\baselineskip,
headfont = \normalfont\bfseries,
notefont = \mdseries, 
notebraces = {(}{)},
bodyfont = \normalfont,
postheadspace = .5em,
headpunct = .,
qed = \qedsymbol
]{solutionstyle}

\renewcommand*{\backref}[1]{}
\renewcommand*{\backrefalt}[4]{
	\ifcase #1 Not cited.%
	\or $\specialuparrow$#2%
	\else $\specialuparrow$#2%
	\fi%
}

\declaretheorem[name = Conjecture, numberwithin = section, style = bolditalic, refname = {Conjecture,Conjectures}, Refname = {Conjecture,Conjectures}]{conjecture}

\declaretheorem[name = Lemma, numberlike = conjecture, style = bolditalic, refname = {Lemma,Lemmas}, Refname = {Lemma,Lemmas}]{lemma}

\declaretheorem[name = Theorem, numberlike = conjecture, style = bolditalic, refname = {Theorem,Theorems}, Refname = {Theorem,Theorems}]{theorem}
\declaretheorem[name = Theorem, numbered = no, style = bolditalic, refname = {Theorem,Theorems}, Refname = {Theorem,Theorems}]{theorem*}
\declaretheorem[name = Question, numberlike = conjecture, style = bolditalic, refname = {Question,Questions}, Refname = {Question,Questions}]{question}

\declaretheorem[name = Remark, numberlike = conjecture, style = boldnormal, refname = {Remark,Remarks}, Refname = {Remark,Remarks}]{remark}

\declaretheorem[name = Proof, numbered = no, style = proofstyle, refname = {Proof,Proofs}, Refname = {Proof,Proofs}]{Proof}

\newlist{property}{enumerate}{1}
\setlist[property]{label = \arabic*., ref = \arabic*}
\crefname{propertyi}{Property}{Properties}

\renewcommand{\epsilon}{\varepsilon}

\DeclarePairedDelimiter{\abs}{\lvert}{\rvert}

\DeclarePairedDelimiter{\ceil}{\lceil}{\rceil}
\DeclarePairedDelimiter{\floor}{\lfloor}{\rfloor}

\DeclarePairedDelimiter{\set}{\{}{\}}

\newcommand*{\cF}{\mathcal{F}}

\newcommand*{\cH}{\mathcal{H}}

\newcommand*{\cO}{\mathcal{O}}
\newcommand*{\cP}{\mathcal{P}}

\newcommand{\defn}[1]{\textcolor{Maroon}{\emph{#1}}}

\begin{document}

\author{Freddie Illingworth\footnotemark[1]}

\title{\bf Minimum degree stability of $H$-free graphs}

\date{12 May 2023}

\maketitle

\begin{abstract}
	Given an $(r + 1)$-chromatic graph $H$, the fundamental edge stability result of Erd\H{o}s and Simonovits says that all $n$-vertex $H$-free graphs have at most $(1 - 1/r + o(1)) \binom{n}{2}$ edges, and any $H$-free graph with that many edges can be made $r$-partite by deleting $o(n^{2})$ edges.
	
	Here we consider a natural variant of this -- the minimum degree stability of $H$-free graphs. In particular, what is the least $c$ such that any $n$-vertex $H$-free graph with minimum degree greater than $cn$ can be made $r$-partite by deleting $o(n^{2})$ edges? We determine this least value for all 3-chromatic $H$ and for very many non-3-colourable $H$ (all those in which one is commonly interested) as well as bounding it for the remainder. This extends the Andr\'{a}sfai-Erd\H{o}s-S\'{o}s theorem and work of Alon and Sudakov.
\end{abstract}

\renewcommand{\thefootnote}{\fnsymbol{footnote}} 


\footnotetext[1]{Mathematical Institute, University of Oxford (\textsf{illingworth@maths.ox.ac.uk}). Research carried out while at DPMMS, University of Cambridge. Research supported by EPSRC grant 2114463.}

\renewcommand{\thefootnote}{\arabic{footnote}} 

\section{Introduction}

A central theme in extremal graph theory is the structure of graphs which do not contain some fixed subgraph $H$. Some of the oldest questions in the area are what is the greatest number of edges such a graph can have, what is the structure of these extremal graphs and how stable is this structure. The classical result is that of Erd\H{o}s and Simonovits.

\begin{theorem}[Erd\H{o}s \& Simonovits,~\cite{Erdos1967, Erdos1968, Simonovits1968}]\label{ES}
	Fix a graph $H$ with chromatic number $r + 1$. If $G$ is an $H$-free graph with $n$ vertices, then $e(G) \leqslant (1 - 1/r + o(1)) \binom{n}{2}$. Furthermore, if $G$ has $(1 - 1/r + o(1)) \binom{n}{2}$ edges, then $G$ can be obtained from the $r$-partite Tur\'{a}n graph, $T_{r}(n)$, by deleting and adding $o(n^{2})$ edges.
\end{theorem}

In particular, all $H$-free graphs with close to the extremal number of edges are close to $(\chi(H) - 1)$-partite. In place of the number of edges, it is natural to consider the structure of $H$-free graphs with large minimum degree and whether they are close to $(\chi(H) - 1)$-partite. In 1974, Andr\'{a}sfai, Erd\H{o}s and S\'{o}s~\cite{AES1974} did this for $H$ being an $(r + 1)$-clique. Their seminal result was the following.

\begin{theorem}[Andr\'{a}sfai-Erd\H{o}s-S\'{o}s,~\cite{AES1974}]\label{AES}
	Let $r \geqslant 2$ and $G$ be a $K_{r + 1}$-free graph with $n$ vertices and minimum degree greater than
	\begin{equation*}
		\biggl(1 - \frac{1}{r - 1/3}\biggr) n.
	\end{equation*}
	Then $G$ is $r$-colourable. Furthermore $1 - 1/(r - 1/3)$ is tight.
\end{theorem}

Note that the $K_{r + 1}$-free graph with greatest minimum degree (and most edges), the Tur\'{a}n graph, $T_{r}(n)$, has minimum degree $(1 - 1/r)n$ so the theorem gives information about graphs whose minimum degrees are $\Omega(n)$ away from that of the extremal graph -- we are far from the regime of \cref{ES}.

We now consider \cref{AES} from the perspective of minimum degree stability. The graph consisting of a 5-cycle joined to a complete $(r - 2)$-partite graph with three vertices in each part shows that $1 - 1/(r - 1/3)$ is tight. Furthermore, $n$-vertex balanced blow-ups (defined in \cref{sec:notation}) of this graph are $K_{r + 1}$-free, have minimum degree $\floor{(1 - 1/(r - 1/3)) n}$, and require the deletion of $\Omega(n^{2})$ edges to be made $r$-partite. This suggests the most basic problem for the minimum degree stability of $H$-free graphs: given a $(r + 1)$-chromatic graph $H$, determine
\begin{equation*}
	\begin{split}
		\delta_{H} = \inf \set{c \colon &\textrm{if }\abs{G} = n, \, \delta(G) \geqslant c n \textrm{ and } G \textrm{ is } H\textrm{-free}, \\
		& \textrm{then } G \textrm{ can be made } r\textrm{-partite by deleting } o(n^{2}) \textrm{ edges}}.
	\end{split}
\end{equation*}
This is the analogue of the edge stability question answered by \cref{ES}. The reason for allowing the omission of $o(n^{2})$ edges is the same for both questions: there are many $H$ for which the $H$-free graph with most edges or greatest minimum degree is not $r$-partite (but of course is very close to $r$-partite -- it is close to the Tur\'{a}n graph). This is explored further in \cref{sec:context}.

It follows from our previous discussion that
\begin{equation*}
	\delta_{K_{r + 1}} = 1 - \frac{1}{r - 1/3}.
\end{equation*}
As pointed out by Alon and Sudakov~\cite{AlonSudakov2006AES}, a standard application of Szemer\'{e}di's regularity lemma~\cite{Szemeredi1978regularity} shows that
\begin{equation*}
	\delta_{K_{r + 1}(t)} = 1 - \frac{1}{r - 1/3},
\end{equation*}
where $K_{r + 1}(t)$ is the complete $(r + 1)$-partite graph with $t$ vertices in each part (so is a Tur\'{a}n graph). In fact, Alon and Sudakov showed that one can do better than deleting $o(n^{2})$ edges.
\begin{theorem}[Alon-Sudakov,~\cite{AlonSudakov2006AES}]\label{AS}
	Let $r \geqslant 2$ and $t \geqslant 1$ be integers, let $\varepsilon > 0$, and set $\rho = 1/(4r^{2/3}t)$. The following holds for all sufficiently large $n$. If $G$ is a $K_{r + 1}(t)$-free graph on $n$ vertices and with minimum degree at least $(1 - 1/(r - 1/3) + \varepsilon)n$, then one can delete $\cO_{r, t, \varepsilon}(n^{2 - \rho})$ edges to make $G$ $r$-colourable.
\end{theorem}
Subsequently, Allen~\cite{Allen2010} found a more direct proof (with no use of the regularity lemma) that yields optimal (to within a constant factor) estimates. Given the central role played by $K_{r + 1}(t)$ in edge stability, it might be tempting to believe that this should determine $\delta_{H}$ for general non-bipartite $H$. Of course, if $H$ is a graph with chromatic number $\chi(H) = r + 1 \geqslant 3$, then $H$ is a subgraph of $K_{r + 1}(t)$ for some $t$ and so any $H$-free graph is $K_{r + 1}(t)$-free. Thus,
\begin{equation}
	\delta_{H} \leqslant \delta_{K_{r + 1}(t)} = 1 - \frac{1}{r - 1/3}.\label{eq:AS}
\end{equation}
However, the inequality may be strict. In the case of edge stability, the Tur\'{a}n graph $T_{r}(n)$ has $(1 - 1/r) \binom{n}{2}$ edges and does not contain $H$, as $H$ has chromatic number $r + 1$. Here, however, it is blow-ups of $K_{r - 2} + C_{5}$ (mentioned above) which show that $\delta_{K_{r + 1}(t)} \geqslant 1 - 1/(r - 1/3)$. These are $K_{r + 1}$-free, but need not be $H$-free -- for a simple example, consider $r = 2$ and $H$ an odd cycle. This observation highlights the following important notion. For two graphs $H$ and $F$ we say there is a \defn{homomorphism $H \to F$}, writing $H \to F$ for short, if $H$ is a subgraph of a blow-up of $F$ -- this is discussed more comprehensively in \cref{sec:notation}. Minimum degree stability is nuanced: $\delta_{H}$ is determined not just by the chromatic number of $H$ but also by its finer structural properties. For 3-chromatic graphs, the situation is fairly straightforward.

\begin{restatable}[$\delta_{H}$ for 3-chromatic $H$]{theorem}{deltaHthree}\label{deltaH3}
	Let $H$ be a 3-chromatic graph. There is a smallest positive integer $g$ for which there is \textbf{no} homomorphism $H \to C_{2g + 1}$. Then
	\begin{equation*}
		\delta_{H} = \frac{2}{2g + 1}.
	\end{equation*}
\end{restatable}

Thus, for 3-chromatic $H$, $\delta_{H}$ is determined by the first odd cycle to which there is not a homomorphism from $H$. Next we turn to graphs which are not 3-colourable. We will determine $\delta_{H}$ for very many $H$ (indeed, all those in which one is commonly interested) and bound it for the remainder.
\begin{theorem}[$\delta_{H}$ for $H$ not 3-colourable]\label{deltaH4}
	There is a sequence of eleven graphs $(F_{g})_{1 \leqslant g \leqslant 11}$ \textnormal{(}described explicitly in \cref{sec:graphs}\textnormal{)} and constants $(c_{g})_{1 \leqslant g \leqslant 11}$ such that the following holds. Let $r \geqslant 3$ be an integer. If $H$ is an $(r + 1)$-chromatic graph and $g$ is minimal such that there is \textbf{no} homomorphism $H \to K_{r - 3} + F_{g}$, then
	\begin{equation*}
		\delta_{H} = 1 - \frac{1}{r - 1 + c_{g}}.
	\end{equation*}
	If there is a homomorphism from $H$ to each of the eleven $K_{r - 3} + F_{g}$, then there is a least $g$ for which there is \textbf{no} homomorphism $H \to K_{r - 2} + C_{2g + 1}$\textnormal{:} $\delta_{H}$ satisfies the bounds
	\begin{equation*}
		1 - \frac{1}{r - 1} < 1 - \frac{1}{r - 1 + 2/(2g - 1)} \leqslant \delta_{H} \leqslant 1 - \frac{1}{r - 1 + 1/7}.
	\end{equation*}
\end{theorem}

Again, there is a sequence of graphs such that $\delta_{H}$ is determined by the first one to which there is no homomorphism from $H$. We will define the $F_{g}$ and $c_{g}$ explicitly in \cref{sec:graphs} -- we mention for now that the $F_{g}$ are all 4-chromatic graphs on few vertices.

The only $(r + 1)$-chromatic graphs whose $\delta_{H}$ value is not determined by \cref{deltaH4} are those which have a homomorphism to each of the eleven $K_{r - 3} + F_{g}$. Such a graph would be rare.

\subsection{Wider context -- chromatic profiles}\label{sec:context}

We give a second view of $\delta_{H}$, placing it within a whole spectrum of structural constants relating to $H$-free graphs. The \defn{chromatic profile of a graph $H$} is the sequence of values for $k = 1, 2, \dotsc$ of
\begin{equation*}
	\delta_{\chi}(H, k) = \inf\set{c \colon \text{if }\delta(G) \geqslant c \abs{G} \text{ and } G \text{ is } H\text{-free}, \text{then } G \text{ is $k$-colourable}}.
\end{equation*}
The question of determining this was first asked by Erd\H{o}s and Simonovits~\cite{ErdosSimonovits1973} in 1973. Note that the Andr\'{a}sfai-Erd\H{o}s-S\'{o}s theorem says that $\delta_{\chi}(K_{r + 1}, r) = 1 - 1/(r - 1/3)$. The chromatic profile of triangles has been studied extensively~\cite{AES1974,Brandt1999,CJK1997,Haggkvist1982,Jin1995,Luczak2006,Thomassen2002} and was finally settled by Brandt and Thomass\'{e}~\cite{BrandtThomasse2005}. These results were extended to general cliques~\cite{GoddardLyle2010,Nikiforov2010}.

Much less is known about the chromatic profile of non-complete graphs with Erd\H{o}s and Simonovits describing the behaviour as ``too complicated''. In part this is because the exact structure and chromatic number of the $n$-vertex $H$-free graphs with highest minimum degree or most edges is unknown. Moreover, degenerate examples abound. Consider $H = K_{3}(2)$ and let $G$ be the complete bipartite Tur\'{a}n graph $T_2(n)$ with a graph $F$ inserted into one of the parts. Now if $F$ does not contain a 4-cycle, then $G$ is $H$-free. Furthermore, $\chi(G) \geqslant \chi(F)$ and $\delta(G) \geqslant \floor{n/2}$. Thus, taking $F$ to have girth at least five and arbitrarily large chromatic number shows that $\delta_{\chi}(H, k) \geqslant 1/2$ for all $k$. But, of course, all $H$-free graphs have at most $(1/2 + o(1)) \binom{n}{2}$ edges and so minimum degree at most $(1/2 + o(1))n$. In particular, $1/2 = \delta_{\chi}(H, k)$ for all $k$. This is unsatisfying, failing to capture the macroscopic behaviour of $H$-free graphs with large minimum degree. All these graphs are close to (within $o(n^2)$ edge of) being bipartite and their high chromatic number is rather artificial. For these reasons, another natural notion of the structure of $H$-free graphs with large minimum degree is the \defn{approximate chromatic profile}. This is
\begin{equation*}
	\begin{split}
		\delta^{\ast}_{\chi}(H, k) = \inf \set{c \colon &\textrm{if }\abs{G} = n, \, \delta(G) \geqslant c n \text{ and } G \textrm{ is } H\textrm{-free}, \\
		& \textrm{then } G \textrm{ can be made } k\textrm{-colourable by deleting } o(n^{2}) \textrm{ edges}}.
	\end{split}
\end{equation*}
Note in passing that $\delta_{\chi}(H, k) \geqslant \delta^{\ast}_{\chi}(H, k)$, although equality need not occur. Indeed, we have just seen that the chromatic profile of $H = K_{3}(2)$ is the constant $1/2$ sequence, while \cref{profiles} below shows that $H$'s approximate chromatic profile is the same as the chromatic profile of triangle-free graphs.

If $H$ is $(r + 1)$-chromatic, then the $r$-partite Tur\'{a}n graph $T_{r}(n)$ is $H$-free and cannot be made $(r - 1)$-colourable without deleting $\Omega(n^{2})$ edges. Also, any $n$-vertex $H$-free graph has at most $(1 - 1/r + o(1))\binom{n}{2}$ edges and so has minimum degree at most $(1 - 1/r + o(1)) n$. Thus, for all $k \leqslant r - 1$ we have $\delta^{\ast}_{\chi}(H, k) = 1 - 1/r$. In particular, the first interesting threshold in the approximate chromatic profile of $H$ is
\begin{equation*}
	\delta^{\ast}_{\chi}(H, \chi(H) - 1)
\end{equation*}
which is exactly $\delta_{H}$.

For a family of graphs $\cF$, one can make the more general definition
\begin{equation*}
	\begin{split}
		\delta^{\ast}_{\chi}(\cF, k) = \inf \set{c \colon &\textrm{if }\abs{G} = n, \, \delta(G) \geqslant c n \text{ and } G \in \cF, \\
		& \textrm{then } G \textrm{ can be made } k\textrm{-colourable by deleting } o(n^{2}) \textrm{ edges}}.
	\end{split}
\end{equation*}
This again satisfies the inequality $\delta_{\chi}(\cF, k) \geqslant \delta^{\ast}_{\chi}(\cF, k)$. There is a natural class of families where equality occurs. We say \defn{$\cF$ is closed under taking blow-ups} if any blow-up of any member of $\cF$ is also in $\cF$. Examples of such families include $K_{t}$-free graphs and $t$-colourable graphs.

\begin{restatable}{theorem}{profiles}\label{profiles}
	Let $\cF$ be a family of graphs that is closed under taking blow-ups. For any positive integer $k$,
	\begin{equation*}
		\delta_{\chi}^{\ast}(\cF, k) = \delta_{\chi}(\cF, k).
	\end{equation*}
\end{restatable}

For non-complete $H$, the family of $H$-free graphs is not closed under taking blow-ups. However, there is a natural family of graphs which is closed under taking blow-ups and whose chromatic profile is the same as the approximate chromatic profile of $H$. Define \defn{$H$-hom} to be the family of \defn{$H$-homomorphism-free graphs}, that is, those graphs to which there is no homomorphism from $H$. Then
\begin{equation*}
	\begin{split}
		\delta_{\chi}(H\text{-hom}, k) = \inf\set{c \colon &\text{if }\delta(G) \geqslant c \abs{G} \text{ and there is no homomorphism } H \to G, \\
		& \text{then } G \text{ is $k$-colourable}},
	\end{split}
\end{equation*}
is the chromatic profile of this family. In \cref{sec:approx}, we will show that this is identical to the approximate chromatic profile of $H$ and so $\delta_{H} = \delta_{\chi}(H\text{-hom}, \chi(H) - 1)$.

\begin{restatable}{theorem}{profileshom}\label{profileshom}
	For any graph $H$ and any positive integer $k$,
	\begin{equation*}
		\delta^{\ast}_{\chi}(H, k) = \delta_{\chi}(H\text{\emph{-hom}}, k).
	\end{equation*}
\end{restatable}

\subsection{The graphs in \texorpdfstring{\cref{deltaH4}}{Theorem 4} and some motivation}\label{sec:graphs}

Various graphs appear as $F_{g}$ in the statement of \cref{deltaH4}. Here we define them explicitly and provide some motivation for their presence. The $F_{g}$ and $c_{g}$ are given in \cref{table:graphs} and the $F_{g}$ are shown in \cref{fig:graphs}.

\begin{table}[ht]
	\centering
	\begin{tabular}{c c c c c c c c c c c c }
		\toprule
		$g$ & 1 & 2 & 3 & 4 & 5 & 6 & 7 & 8 & 9 & 10 & 11 \\
		\midrule
		$F_{g}$  & $W_{5}$ & $W_{7}$ & $\overline{C}_{7}$ & $W_{9}$ & $H_{2}^{+}$ & $W_{11}$ & $H_{2}$ & $W_{13}$ & $T_{0}$ & $W_{15}$ & $H_{1}^{++}$ \\ \addlinespace
		$c_{g}$ & $\frac{2}{3}$ & $\frac{2}{5}$ & $\frac{1}{3}$ & $\frac{2}{7}$ & $\frac{1}{4}$ & $\frac{2}{9}$ & $\frac{1}{5}$ & $\frac{2}{11}$ & $\frac{1}{6}$ & $\frac{2}{13}$ & $\frac{1}{7}$ \\
		\bottomrule
	\end{tabular}
	\caption{$F_{g}$ and $c_{g}$}\label{table:graphs}
\end{table}

The graph $W_{k}$ (called a $k$-wheel) is a single vertex joined to a $k$-cycle. The graph $\overline{C}_{7}$ is the complement (and also the square) of a 7-cycle. The graph $H_{2}$ is obtained from $\overline{C}_{7}$ by deleting an edge (while maintaining 4-chromaticness) and $H_{2}^{+}$ is obtained from $H_{2}$ by adding a vertex of degree three. The graph $H_{1}^{++}$ can be obtained from $H_{2}$ by deleting an edge and adding two degree three vertices. Finally, $T_{0}$ is a 7-cycle (the outer cycle) together with two vertices each joined to six of the seven vertices in the outer cycle (with the `seventh' vertices two apart) as well as a vertex of degree three.

\begin{figure}[H]
	\centering
	\begin{subfigure}{.24\textwidth}
		\centering
		\begin{tikzpicture}
			\foreach \pt in {0,1,...,6} 
			{
				\tkzDefPoint(\pt*360/7 + 90:1){v_\pt}
			} 
			\tkzDefPoint(0,0){u}
			\tkzDrawPolySeg(v_0,v_1,v_2,v_3,v_4,v_5,v_6,v_0) 
			\foreach \pt in {0,1,...,6} 
			{
				\tkzDrawSegment(u,v_\pt)
			} 
			\tkzDrawPoints(v_0,v_...,v_6)
			\tkzDrawPoint(u)
		\end{tikzpicture}
		\subcaption*{$W_{7}$}
	\end{subfigure}
	\begin{subfigure}{.24\textwidth}
		\centering
		\begin{tikzpicture}
			\foreach \pt in {0,1,...,6} 
			{
				\tkzDefPoint(\pt*360/7 + 90:1){v_\pt}
			} 
			\tkzDrawPolySeg(v_0,v_1,v_2,v_3,v_4,v_5,v_6,v_0) 
			\tkzDrawPolySeg(v_0,v_2,v_4,v_6,v_1,v_3,v_5,v_0)
			\tkzDrawPoints(v_0,v_...,v_6)
		\end{tikzpicture}
		\subcaption*{$C_{7}^{2} = \overline{C}_{7}$}
	\end{subfigure}
	\begin{subfigure}{.24\textwidth}
		\centering
		\begin{tikzpicture}
			\foreach \pt in {0,1,...,6} 
			{
				\tkzDefPoint(\pt*360/7 + 90:1){v_\pt}
			} 
			\tkzDefPoint(0,0){u}
			\tkzDrawPolySeg(v_0,v_1,v_2,v_3,v_4,v_5,v_6,v_0)
			\tkzDrawPolySeg(v_1,v_3,v_5,v_0, v_2,v_4,v_6)
			\tkzDrawSegments(u,v_0 u,v_2 u,v_5)
			\tkzDrawPoints(v_0,v_...,v_6)
			\tkzDrawPoint(u)
		\end{tikzpicture}
		\subcaption*{$H_{2}^{+}$}
	\end{subfigure}
	\begin{subfigure}{.24\textwidth}
		\centering
		\begin{tikzpicture}
			\foreach \pt in {0,1,...,6} 
			{
				\tkzDefPoint(\pt*360/7 + 90:1){v_\pt}
			} 
			\tkzDrawPolySeg(v_0,v_1,v_2,v_3,v_4,v_5,v_6,v_0)
			\tkzDrawPolySeg(v_1,v_3,v_5,v_0, v_2,v_4,v_6)
			\tkzDrawPoints(v_0,v_...,v_6)
		\end{tikzpicture}
		\subcaption*{$H_{2}$}
	\end{subfigure}
	
	\bigskip
	
	\begin{subfigure}{.24\textwidth}
		\centering
		\begin{tikzpicture}
			\foreach \pt in {0,1,...,6} 
			{
				\tkzDefPoint(\pt*360/7 + 90:1){v_\pt}
			}
			\tkzDefPoint(0,0.33){t}
			\tkzDefPoint(0.46,-0.4){u_1}
			\tkzDefPoint(-0.46,-0.4){u_6}
			\tkzDrawPolySeg(v_0,v_1,v_2,v_3,v_4,v_5,v_6,v_0) 
			\tkzDrawSegments(t,v_0 t,u_1 t,u_6)
			\foreach \pt in {0,2,3,4,5,6}
			{
				\tkzDrawSegment(u_1,v_\pt)
			}
			\foreach \pt in {0,1,2,3,4,5}
			{
				\tkzDrawSegment(u_6,v_\pt)
			}
			\tkzDrawPoints(v_0,v_...,v_6)
			\tkzDrawPoints(t,u_1,u_6)
		\end{tikzpicture}
		\subcaption*{$T_{0}$}
	\end{subfigure}
	\begin{subfigure}{.24\textwidth}
		\centering
		\begin{tikzpicture}
			\foreach \pt in {0,1,...,6} 
			{
				\tkzDefPoint(\pt*360/7 + 90:1){v_\pt}
			}
			\tkzDefPoint(-0.33,0){ul}
			\tkzDefPoint(0.33,0){ur}
			
			\tkzDrawPolySeg(v_0,v_1,v_2,v_3,v_4,v_5,v_6,v_0) 
			\tkzDrawPolySeg(v_3,v_5,v_0,v_2,v_4)
			\tkzDrawPolySeg(v_6,v_1)
			\tkzDrawSegments(ul,v_0 ul,v_2 ul,v_3 ur,v_0 ur,v_5 ur,v_3)
			\tkzDrawPoints(v_0,v_...,v_6)
			\tkzDrawPoints(ul,ur)
		\end{tikzpicture}
		\subcaption*{$H_{1}^{++}$}
	\end{subfigure}
	\caption{$F_{g}$}\label{fig:graphs}
\end{figure}

The sequence $F_{g}$ is slightly unusual. Firstly the graphs do not increase in size. Secondly it is not always true that there is a homomorphism $F_{g + 1} \to F_{g}$ and so it is, for example, possible for a graph to have a homomorphism to both $F_{3}$ and $F_{5}$ but not to $F_{4}$.

We now motivate why it is these graphs that are the $F_{g}$. The intuitive explanation for \cref{deltaH3} is that the main obstacle for being close to (that is, within $o(n^{2})$ edges of) bipartite is containing some blow-up (defined in \cref{sec:notation}) of an odd cycle. That odd cycle must be consistent with being $H$-free (in particular, the blow-up of the odd cycle must be $H$-free) and hence it is the first odd cycle to which there is no homomorphism from $H$ that determines $\delta_{H}$. 

Now consider the $r = 3$ version of \cref{deltaH4}: we are interested in which graphs' blowups are the main obstacles for being close to tripartite. Given the importance of odd cycles for being far from bipartite it seems natural that odd wheels would be obstacles here and indeed six of the $F_{g}$ are odd wheels. The other five $F_{g}$ do not contain any odd wheels and so all their neighbourhoods are bipartite -- they are \defn{locally bipartite} graphs. These observations suggest we should pay attention to 4-chromatic locally bipartite graphs as these may be obstacles for being close to tripartite. More generally, call a graph \defn{$a$-locally bipartite} if the common neighbourhood of every $a$-clique is bipartite (so that 1-locally bipartite graphs are exactly locally bipartite). The following result, which appeared in~\cite{Illingworth2022localbpart}, provides $a$-locally bipartite graphs which are obstacles for being $(a + 2)$-colourable. This will be used as part of our proof of \cref{deltaH4}.

\begin{theorem}[$a$-locally bipartite graphs]\label{spec4alocalbip}
	Let $G$ be an $a$-locally bipartite graph.
	\begin{itemize}[noitemsep]
		\item If $\delta(G) > (1 - 1/(a + 4/3)) \cdot \abs{G}$, then $G$ is $(a + 2)$-colourable.
		\item If $\delta(G) > (1 - 1/(a + 5/4)) \cdot \abs{G}$, then $G$ is either $(a + 2)$-colourable or contains $K_{a - 1} + \overline{C}_{7}$.
		\item If $\delta(G) > (1 - 1/(a + 6/5)) \cdot \abs{G}$, then $G$ is either $(a + 2)$-colourable or contains $K_{a - 1} + \overline{C}_{7}$ or $K_{a - 1} + H_{2}^{+}$.
		\item If $\delta(G) > (1 - 1/(a + 7/6)) \cdot \abs{G}$, then $G$ is either $(a + 2)$-colourable or contains $K_{a - 1} + H_{2}$.
		\item If $\delta(G) > (1 - 1/(a + 8/7)) \cdot \abs{G}$, then $G$ is either $(a + 2)$-colourable or contains $K_{a - 1} + H_{2}$ or $K_{a - 1} + T_{0}$.
	\end{itemize}
\end{theorem}

One might ask whether there are other sequences $F'_{g}$ and $c'_{g}$ for which \cref{deltaH4} holds. The fact that $\delta_{H}$ is a fixed number has two corollaries. Firstly, it must be the case that $c'_{g} = c_{g}$ for $1 \leqslant g \leqslant 11$. Secondly, any graph with a homomorphism to each of $F_{1}, F_{2}, \dotsc, F_{g}$ must also have a homomorphism to each of $F_{1}, F_{2}, \dotsc, F_{g - 1}, F'_{g}$ and vice versa. It seems likely that the $F_{g}$ are the minimal graphs satisfying this and so form the ``canonical'' sequence, but proving this is not straightforward.

\subsection{Notation}\label{sec:notation}

Given a graph $G$, a \defn{blow-up} of $G$, is a graph obtained by replacing each vertex $v$ of $G$ by a non-empty independent set $I_{v}$ and each edge $uv$ by a complete bipartite graph between classes $I_{u}$ and $I_{v}$. We say we have \defn{blown-up a vertex $v$ by $n$} if $\abs{I_v} = n$. It is often helpful to think of this as weighting vertex $v$ by $n$.

A blow-up is \defn{balanced} if the independent sets $(I_{v})_{v \in G}$ are as equal in size as possible. We use $G(t)$ to denote the graph obtained by blowing-up each vertex of $G$ by $t$, i.e.\ $G(t)$ is the balanced blow-up of $G$ on $t \abs{G}$ vertices. For example, the balanced blow-up of the $r$-clique, $K_{r}$, on $n$ vertices is exactly the Tur\'{a}n graph, $T_{r}(n)$. We note in passing that a graph has the same chromatic and clique number as any of its blow-ups. Furthermore, if $H$ is a blow-up of $G$, then $G$ is $a$-locally bipartite if and only if $H$ is.

Given two graphs $G$ and $H$, the \defn{join} of $G$ and $H$, denoted by \defn{$G + H$}, is the graph obtained by taking disjoint copies of $G$ and $H$ and joining each vertex of the copy of $G$ to each vertex of the copy of $H$. Note that the chromatic and clique numbers of $G + H$ are the sum of the chromatic and clique numbers of $G$ and $H$.

There is a homomorphism from a graph $G$ to a graph $H$, written \defn{$G \to H$}, if there is a map $\varphi \colon V(G) \to V(H)$ such that for any edge $uv$ of $G$, $\varphi(u)\varphi(v)$ is an edge of $H$. There is a homomorphism $G \to H$ if and only if $G$ is a subgraph of some blow-up of $H$. In particular, if $G \to H$, then $\chi(G) \leqslant \chi(H)$ and moreover if $H$ is $a$-locally bipartite, then $G$ is also.

\subsection{Tools}\label{sec:tools}

We will make great use of Szemer\'{e}di's regularity lemma~\cite{Szemeredi1978regularity} together with some associated machinery which we describe here. Let $(X, Y)$ be a pair of vertex subsets of graph $G$. We use $d(X, Y) = e(X, Y) \abs{X}^{-1} \abs{Y}^{-1}$ to denote the \defn{density} between $X$ and $Y$. The pair $(X, Y)$ is \emph{$\varepsilon$-regular} if 
\begin{equation*}
	\abs{d(U, V) - d(X, Y)} \leqslant \varepsilon
\end{equation*}
for all $U \subset X$, $V \subset Y$ with $\abs{U} \geqslant \varepsilon \abs{X}$ and $\abs{V} \geqslant \varepsilon \abs{Y}$. A partition $\cP = V_{0} \cup V_{1} \cup \dotsb \cup V_{k}$ of $V(G)$ is an \defn{$\varepsilon$-regular partition} if:
\begin{itemize}[noitemsep]
	\item $V_{1}$, $V_{2}$, \ldots, $V_{k}$ all have equal size and $\abs{V_0} \leqslant \varepsilon \abs{G}$,
	\item for all but at most $\varepsilon \binom{k}{2}$ pairs $ij$ ($i < j$), the pair $(V_{i}, V_{j})$ is $\varepsilon$-regular. 
\end{itemize}
Szemer\'{e}di's celebrated result is that, for every positive integer $\ell$ and $\varepsilon > 0$, there is some $L = L(\ell, \varepsilon)$ such that every graph with at least $\ell$ vertices has an $\varepsilon$-regular partition into at least $\ell$ but at most $L$ parts. We will need a version of Szemer\'{e}di's regularity lemma which works well with minimum degrees.

Fix $\varepsilon > 0$ and some $\lambda \geqslant 0$. Suppose we have a graph $G$ with some $\varepsilon$-regular partition $\cP = V_{0} \cup V_{1} \cup \dotsb \cup V_{k}$. These induce what is called a \defn{reduced graph} $R(\cP, \varepsilon, \lambda)$: this has vertex set $\set{1, 2, \dotsc, k}$ with vertex $i$ joined to vertex $j$ exactly if the pair $(V_{i}, V_{j})$ is $\varepsilon$-regular and $d(V_{i}, V_{j}) \geqslant \lambda$. Note when $\lambda = 0$ this graph will have at least $(1 - \varepsilon) \binom{k}{2}$ edges by the definition of an $\varepsilon$-regular partition. We will make use of the following version of the regularity lemma which is an immediate corollary of Theorem~1.10 in Koml\'{o}s and Simonovits's survey of the subject~\cite{KomlosSimonovits1996}.

\begin{lemma}[Szemer\'{e}di's regularity lemma, minimum degree form]\label{SRL}
	Let $\varepsilon > 0$, $\lambda, \delta \in [0, 1]$ and $\ell$ be a positive integer. There is a positive integer $L$ such that the following holds for all $n \geqslant \ell$. If $G$ is graph on $n$ vertices with $\delta(G) \geqslant \delta n$, then $G$ has some $\varepsilon$-regular partition $\cP = V_{0} \cup V_{1} \cup \dotsb \cup V_{k}$ with $k$ between $\ell$ and $L$ such that the corresponding reduced graph $R(\cP, \varepsilon, \lambda)$ has minimum degree at least $(\delta - \varepsilon - \lambda) k$.
\end{lemma}
The point of the reduced graph is that if it contains some structure, then we can find a large structure in $G$, by using a building lemma (for example, see Theorem 2.1 in~\cite{KomlosSimonovits1996}).
\begin{lemma}[graph-building lemma]\label{Szembuilding}
	Let $H$ be a graph \textup{(}on vertex set $\set{1, 2, \dotsc, \abs{H}}$\textup{)}, $t$ a positive integer and $\lambda \in (0, 1)$. For all sufficiently small $\varepsilon > 0$ the following holds. Suppose $V_{1}$, \ldots, $V_{\abs{H}}$ are sufficiently large pairwise disjoint vertex sets with $(V_{i}, V_{j})$ $\varepsilon$-regular of density at least $\lambda$ for each $ij \in E(H)$. Then we can find a copy of $H(t)$ with each blown-up vertex in the corresponding $V_{i}$.
\end{lemma}

These two lemmas work well together: given a large graph $G$ with minimum degree $\delta \abs{G}$, \cref{SRL} shows that $G$ has a corresponding reduced graph $R$ of bounded size and with minimum degree almost $\delta \abs{R}$. If $G$ is $H(t)$-free, then, by \cref{Szembuilding}, $R$ is $H$-free. This may give some structural information about $R$ (e.g.\ it is $r$-colourable) which we will then pull back to $G$.

Finally, odd cycles will play an important role in determining $\delta_{H}$, so we note the following fact about homomorphisms to odd cycles.

\begin{lemma}\label{deltaHlemma4girth}
	Let $g$ be a positive integer. If there is a homomorphism $G \to C_{2g + 1}$, then $G$ contains no odd cycles of length less than $2g + 1$.
\end{lemma}

\begin{Proof}
	Let $\varphi$ be a homomorphism from $G$ to $C_{2g + 1}$. Let $C$ be an odd cycle of $G$. The restriction of $\varphi$ to $C$ gives a homomorphism from $C$ to $\varphi(C)$, so $\chi(\varphi(C)) \geqslant \chi(C) = 3$. Thus, $\varphi(C)$ is the whole of $C_{2g + 1}$, so $\abs{C} \geqslant \abs{\varphi(C)} = 2g + 1$.
\end{Proof}

\section{\texorpdfstring{$\delta_{H}$}{deltaH} for 3-chromatic \texorpdfstring{$H$}{H} -- proof of \texorpdfstring{\cref{deltaH3}}{Theorem 3}}\label{sec:deltaH3}

Let $H$ be a graph with chromatic number three, so there is a homomorphism $H \to K_{3} \cong C_{3}$. Furthermore, $H$ is not bipartite, so contains at least one odd cycle. This, coupled with \cref{deltaHlemma4girth}, means that there is a homomorphism from $H$ to only finitely many odd cycles. Let $g$ be the smallest positive integer for which there is no homomorphism $H \to C_{2g + 1}$. A balanced blow-up of $C_{2g + 1}$ on $n$ vertices is $H$-free and has minimum degree at least $2 \floor{n/(2g + 1)}$. We claim that to make this balanced blow-up bipartite requires the deletion of at least $\floor{n/(2g + 1)}^{2} = \Omega(n^{2})$ edges. Indeed, let the sizes of the $(2g + 1)$ parts in the blow-up be $\ceil{n/(2g + 1)} = x_{1} \geqslant x_{2} \geqslant \dotsb \geqslant x_{2g + 1} = \floor{n/(2g + 1)}$. Then the number of copies of $C_{2g + 1}$ in the blow-up is $x_{1}x_{2} \dotsm x_{2g + 1}$ and each edge lies in at most $x_{1} x_{2} \dotsm x_{2g - 1}$ copies of $C_{2g + 1}$, so to make the blow-up bipartite requires the deletion of at least $x_{2g} x_{2g + 1} \geqslant \floor{n/(2g + 1)}^{2}$ edges. In particular,
\begin{equation*}
	\delta_{H} \geqslant \frac{2}{2g + 1}.
\end{equation*}
We claim that in fact there is equality. Before proving this, we need the following result for odd cycles, which was noted by Andr\'{a}sfai, Erd\H{o}s and S\'{o}s~\cite{AES1974}. For completeness we give a proof.

\begin{lemma}\label{deltaHlemma4oddcycles}
	Let $g \geqslant 2$ be a positive integer. Suppose $G$ is a non-bipartite graph with
	\begin{equation*}
		\delta(G) > \frac{2}{2g + 1} \cdot \abs{G}.
	\end{equation*}
	Then $G$ contains an odd cycle of length less than $2g + 1$.
\end{lemma}

\begin{Proof}
	Let $C$ be the shortest odd cycle in $G$. By minimality, $C$ is induced and no vertex is adjacent to three vertices in $C$. Thus,
	\begin{equation*}
		\abs{C} \cdot \delta(G) \leqslant e(C, G) \leqslant 2 \abs{G},
	\end{equation*}
	and so $\abs{C} \leqslant 2 \abs{G}/\delta(G) < 2g + 1$.
\end{Proof}
We are now ready to prove \cref{deltaH3} which determines $\delta_{H}$ for 3-chromatic $H$.

\deltaHthree*

\begin{Proof}
	The graph $H$ is not bipartite so contains an odd cycle. In particular, if $g$ is such that $2g + 1$ is greater than the length of the shortest odd cycle of $H$, then, by \cref{deltaHlemma4girth}, there is no homomorphism $H \to C_{2g + 1}$. Take a minimal such $g$.
	
	By the opening remarks of this section, a balanced blow-up of $C_{2g + 1}$ on $n$ vertices is $H$-free, has minimum degree at least $2\floor{n/(2g + 1)}$ and requires at least $\Omega(n^{2})$ edges to be deleted to be made bipartite so $\delta_{H} \geqslant 2/(2g + 1)$.
	
	We are left to show that, for all $\eta > 0$, any $n$-vertex $H$-free graph with minimum degree at least
	\begin{equation*}
		\biggl(\frac{2}{2g + 1} + \eta\biggr)n,
	\end{equation*}
	can be made bipartite by deleting at most $\eta n^{2}$ edges when $n$ is sufficiently large. Firstly, note that there is a homomorphism from $H$ to each of $C_{3}$, $C_{5}$, \ldots, $C_{2g - 1}$, so there exists some positive integer $t$ such that $H$ is a subgraph of all of $C_{3}(t)$, $C_{5}(t)$, \ldots, $C_{2g - 1}(t)$.
	
	Fix $n$ and $\ell$ large (chosen later), let $\lambda = \eta /2$ and $\varepsilon > 0$ be sufficiently small. Let $G$ be a $n$-vertex graph with minimum degree at least $(2/(2g + 1) + \eta)n$. By \cref{SRL}, $G$ has some $\varepsilon$-regular partition $\cP = V_{0} \cup V_{1} \cup \dotsb \cup V_{k}$ with $k$ between $\ell$ and $L$ (a constant not depending on $G$ or $n$) such that the reduced graph $R = R(\cP, \varepsilon, \lambda)$ has minimum degree greater than $2k/(2g + 1)$. By \cref{deltaHlemma4oddcycles}, $R$ is either bipartite or contains one of $C_{3}$, $C_{5}$, \ldots, $C_{2g - 1}$.
	
	Applying \cref{Szembuilding}, provided $\varepsilon$ was chosen small enough (in terms of $\lambda$) and $n$ is large enough, either $R$ is bipartite or $G$ contains one of $C_{3}(t)$, $C_{5}(t)$, \ldots, $C_{2g - 1}(t)$. The latter contradicts $G$ being $H$-free and so $R$ is bipartite.
	
	Now, consider deleting from $G$ all edges within each $V_{i}$, the edges incident to $V_{0}$ and all edges between $V_{i}$ and $V_{j}$ when $ij \not\in E(R)$. The resulting graph is a blow-up of $R$, so is bipartite. This process deletes at most
	\begin{align*}
		& k \cdot (n/k)^{2} + \varepsilon n \cdot n + \varepsilon \tbinom{k}{2} (n/k)^{2} + \lambda n^{2} \\
		\leqslant & \ n^{2}(\lambda + 2\varepsilon + 1/k) \\
		\leqslant & \ n^{2}(\eta/2 + 2 \varepsilon + 1/\ell) \leqslant \eta n^{2}
	\end{align*}
	edges, provided that $\ell$ is large enough and $\varepsilon$ is small enough.
\end{Proof}

\begin{remark}
	By \cref{deltaHlemma4girth}, the odd girth of $G$ is at least $2g - 1$ (where $g$ is as in the theorem statement). However, we may not have equality. For example, the Petersen graph has odd girth 5 but has no homomorphism to $C_{5}$, so $\delta_{\text{Petersen}}$ is $2/5$ and not $2/7$.
\end{remark}

\section{Properties of \texorpdfstring{$\delta_{H}$}{deltaH}}\label{sec:deltaHprop}

We take a moment to crystallise the key ingredients of the proof of \cref{deltaH3} (and, in particular, what role the odd cycles played). Fix $r \geqslant 2$ and suppose we have a sequence of graphs $K_{r + 1} = L_{0}, L_{1}, L_{2}, \dotsc, L_{m}$ and a sequence of constants $k_{1}, k_{2}, \dotsc, k_{m}$ where $m$ may be infinity (if both sequences are infinite). The relevant properties this pair of sequences might satisfy are the following.
\begin{property}[noitemsep]
	\item \label{property:1} None of $L_{1}$, $L_{2}$, $L_{3}$, \ldots \ is $r$-colourable.
	\item \label{property:2} No $(r + 1)$-chromatic graph has a homomorphism to all of $L_{1}$, $L_{2}$, \ldots\,.
	\item \label{property:3} For each $g$, if $G$ is an $n$-vertex graph with minimum degree greater than $k_{g} n$, then $G$ is either $r$-colourable or contains at least one of $L_{0}$, $L_{1}$, \ldots, $L_{g - 1}$.
	\item \label{property:4} For each $g$ and any $c < k_{g}$, there is some blow-up $L'_g$ of $L_{g}$ satisfying $\delta(L'_g) \geqslant c \cdot \abs{L'_g}$.
\end{property}
\Cref{deltaH3} corresponds to the sequences $L_{g} = C_{2g + 3}$ and $k_{g} = 2/(2g + 3)$ satisfying all the properties for $r = 2$. As odd cycles are regular, we are able to take $L'_g = L_{g}$ in this case. However, the $L_{g}$ we use later will often be non-regular and furthermore they may have no blow-ups with $\delta(L'_g) = k_{g} \cdot \abs{L'_g}$, but still satisfy \cref{property:4}. One could weaken \cref{property:3} to minimum degree greater than $(k_{g} + o(1)) n$ (which is all we use in our analysis below) but for our purposes this is unnecessary.

Let $H$ be an $(r + 1)$-chromatic graph. Suppose there is no homomorphism $H \to L_{g}$ and \cref{property:1,property:4} hold. Then, for any $c < k_{g}$, let $G$ be a balanced blow-up of $L'_g$ on $n$ vertices: $H \nrightarrow L_{g}$, so $G$ is $H$-free. Furthermore, $G$ has minimum degree at least $\delta(L'_g) \floor{n/ \abs{L'_g}} \geqslant (c - o(1)) n$. By \cref{property:1}, $L'_g$ is not $r$-colourable and hence to make $G$ $r$-colourable requires the deletion of enough edges so that no copy of $L'_g$ remains -- we will show this requires the deletion of $\Omega(n^{2})$ edges. Let the sizes of the $\abs{L'_g}$ parts of $G$ be $\ceil{n/\abs{L'_g}} \geqslant x_{1} \geqslant x_{2} \geqslant \dotsb \geqslant x_{\abs{L'_g}} \geqslant \floor{n/\abs{L'_g}}$. Now, $G$ contains at least $x_{1}x_{2} \dotsm x_{\abs{L'_g}}$ copies of $L'_g$ in which each vertex is in the corresponding part, and each edge of $G$ is in at most $x_{1}x_{2} \dotsm x_{\abs{L'_g} - 2}$ such copies. Hence to remove all copies of $L'_g$ requires the deletion of at least $x_{\abs{L'_g} - 1} \cdot x_{\abs{L'_g}} \geqslant \floor{n/\abs{L'_g}}^{2} = \Omega(n^{2})$ edges. Thus $\delta_{H} \geqslant c$ and so $\delta_{H} \geqslant k_{g}$.

Suppose there is a homomorphism from $H$ to each of $L_{1}$, $L_{2}$, \ldots, $L_{g - 1}$ and \cref{property:3} holds (note that $H \to L_{0}$ also). Then the same regularity argument as in the proof of \cref{deltaH3} shows that $\delta_{H} \leqslant k_{g}$. We now sketch this argument. Let $\eta > 0$ and take a large graph $G$ with minimum degree at least $(k_{g} + \eta) \abs{G}$. We can use \cref{SRL} to get a corresponding reduced graph $R$ with minimum degree greater than $(k_{g} + \eta')\abs{R}$ (some $\eta' \in (0, \eta)$), which, by \cref{property:3}, is either $r$-colourable or contains one of $L_{0}$, $L_{1}$, \ldots, $L_{g - 1}$. In the latter case, we use \cref{Szembuilding} to get a copy of $H$ in $G$ and in the former case we may delete at most $\eta \abs{G}^{2}$ edges from $G$ to leave an $r$-colourable graph.

The upshot of all this is that if \cref{property:1,property:3,property:4} hold, then either $\delta_{H} \leqslant k_{m}$ (if $H \to L_{1}$, $L_{2}$, \ldots, $L_{m}$) or $\delta_{H} = k_{g}$ where $g$ is minimal with $H \nrightarrow L_{g}$. Of course, if \cref{property:2} also holds, then we can determine $\delta_{H}$ for any $(r + 1)$-chromatic $H$. Even without \cref{property:2}, any sequence does determine $\delta_{H}$ for many $H$ and gives an upper bound for the rest. Thus, we are particularly interested in pairs of sequences satisfying \cref{property:1,property:3,property:4}.

We illustrate these remarks by next proving a weak version of \cref{deltaH4}. Note that, when $r = 2$, $\delta_{H}$ could be arbitrarily close to zero (corresponding to $k_{g} \to 0$ as $g \to \infty$). However, this is not the case for larger $r$.

\begin{theorem}\label{deltaHkg}
	If $H$ is a graph with chromatic number $r + 1 \geqslant 3$, then there is a least $g$ for which there is \textbf{no} homomorphism $H \to K_{r - 2} + C_{2g + 1}$ and furthermore
	\begin{equation*}
		1 - \frac{1}{r - 1} < 1 - \frac{1}{r - 1 + 2/(2g - 1)} \leqslant \delta_{H} \leqslant 1 - \frac{1}{r - 1/3}.
	\end{equation*}
\end{theorem}

\begin{Proof}
	The right-hand inequality is just inequality \eqref{eq:AS}. For the middle inequality define the following sequence of graphs: $L_{g} = K_{r - 2} + C_{2g + 3}$ (so $L_{0} = K_{r + 1}$) and let
	\begin{equation*}
		k_{g} = 1 - \frac{1}{r - 1 + 2/(2g + 1)}.
	\end{equation*}
	From the preceding discussion it suffices to show that these sequences satisfy \cref{property:1,property:2,property:4} (but not necessarily \ref{property:3}). Indeed, if these properties hold, then there is a minimal $g$ such that $H \nrightarrow L_{g - 1}$ and so
	\begin{equation*}
		\delta_{H} \geqslant k_{g - 1} = 1 - \frac{1}{r - 1 + 2/(2g - 1)} > 1 - \frac{1}{r - 1}.
	\end{equation*}
	\Cref{property:1} is immediate: each $L_{g}$ has chromatic number $\chi(K_{r - 2}) + \chi(C_{2g + 3}) = r - 2 + 3 = r + 1$. 
	
	Suppose $H$ is some $(r + 1)$-chromatic graph which has a homomorphism to each of the $L_{g}$ and let these homomorphisms be $\varphi_{g} \colon H \to L_{g} = K_{r - 2} + C_{2g + 3}$. Now $\varphi_{g}^{-1}(K_{r - 2})$ is $(r - 2)$-colourable, so $X_{g} = \varphi_{g}^{-1}(C_{2g + 3})$ is not bipartite (as $H$ is not $r$-colourable). In particular, for each $g$, $H[X_{g}]$ is not bipartite but $H[X_{g}] \to C_{2g + 3}$. Thus, by \cref{deltaHlemma4girth}, $H[X_{g}]$ contains an odd cycle of length at least $2g + 3$. Therefore, $H$ contains odd cycles of arbitrary length, which is absurd. Hence, we have \cref{property:2}.
	
	Finally, for \cref{property:4}, let $L'_g = K_{r - 2}(2g + 1) + C_{2g + 3}$, which is a blow-up of $L_{g}$. The graph $L'_g$ has $(2g + 1)(r - 2) + 2g + 3 = (2g + 1)(r - 1) + 2$ vertices and is $[(2g + 1)(r - 2) + 2]$-regular. In particular,
	\begin{equation*}
		k_{g} \cdot \abs{L'_g} = \biggl(1 - \frac{2g + 1}{(r - 1)(2g + 1) + 2}\biggr) \cdot [(2g + 1)(r - 1) + 2] = \delta(L'_g). \qedhere
	\end{equation*}
\end{Proof}

\section{\texorpdfstring{$\delta_{H}$}{deltaH} for general \texorpdfstring{$H$}{H} -- proof of \texorpdfstring{\cref{deltaH4}}{Theorem 4}}\label{sec:deltaH4}

We are now ready to prove \cref{deltaH4}, which we restate here with the explicit $F_{g}$ and $c_{g}$ for convenience.

\begin{theorem}
	Define the sequence of graphs $F_{g}$ and constants $c_{g}$ $(1 \leqslant g \leqslant 11)$ as follows.
	\begin{table}[H]
		\centering
		\begin{tabular}{c c c c c c c c c c c c }
			\toprule
			$g$ & 1 & 2 & 3 & 4 & 5 & 6 & 7 & 8 & 9 & 10 & 11 \\
			\midrule
			$F_{g}$  & $W_{5}$ & $W_{7}$ & $\overline{C}_{7}$ & $W_{9}$ & $H_{2}^{+}$ & $W_{11}$ & $H_{2}$ & $W_{13}$ & $T_{0}$ & $W_{15}$ & $H_{1}^{++}$ \\ \addlinespace
			$c_{g}$ & $\frac{2}{3}$ & $\frac{2}{5}$ & $\frac{1}{3}$ & $\frac{2}{7}$ & $\frac{1}{4}$ & $\frac{2}{9}$ & $\frac{1}{5}$ & $\frac{2}{11}$ & $\frac{1}{6}$ & $\frac{2}{13}$ & $\frac{1}{7}$ \\
		\bottomrule
		\end{tabular}
	\end{table}
	Let $r \geqslant 3$ be an integer. If $H$ is an $(r + 1)$-chromatic graph and $g$ is minimal such that there is \textbf{no} homomorphism $H \to K_{r - 3} + F_{g}$, then 
	\begin{equation*}
		\delta_{H} = 1 - \frac{1}{r - 1 + c_{g}}.
	\end{equation*}
	If there is a homomorphism from $H$ to each of the eleven $K_{r - 3} + F_{g}$, then there is a least $g$ for which there is \textbf{no} homomorphism $H \to K_{r - 2} + C_{2g + 1}$\textnormal{:} $\delta_{H}$ satisfies the bounds
	\begin{equation*}
		1 - \frac{1}{r - 1} < 1 - \frac{1}{r - 1 + 2/(2g - 1)} \leqslant \delta_{H} \leqslant 1 - \frac{1}{r - 1 + 1/7}.
	\end{equation*}
\end{theorem}

\begin{Proof}
	The lower bound
	\begin{equation*}
		1 - \frac{1}{r - 1 + 2/(2g - 1)} \leqslant \delta_{H}
	\end{equation*}
	in the final part of this theorem follows from \cref{deltaHkg}. The discussion following \cref{property:1,property:2,property:3,property:4} shows that what remains is to check that the sequences
	\begin{equation*}
		L_{g} = K_{r - 3} + F_{g}, \qquad k_{g} = 1 - \frac{1}{r - 1 + c_{g}}, \qquad g = 1, 2, \dotsc, 11
	\end{equation*}
	satisfy \cref{property:1,property:3,property:4}. \Cref{property:1} is immediate: for each $g$, $\chi(L_{g}) = r - 3 + \chi(F_{g}) = r + 1$.
	
	We now consider \cref{property:4}. When $F_{g} = W_{2k + 1}$, $c_{g} = 2/(2k - 1)$. For $L_{g} = K_{r - 3} + W_{2k + 1} = K_{r - 2} + C_{2k + 1}$, we take $L'_g = K_{r - 2}(2k - 1) + C_{2k + 1}$. This is a regular blow-up of $L_{g}$ and, as was shown in the proof of \cref{deltaHkg}, satisfies
	\begin{equation*}
		\biggl(1 - \frac{1}{r - 1 + c_{g}}\biggr) \cdot \abs{L'_g} = \delta(L'_g),
	\end{equation*}
	which is exactly $k_{g} \cdot \abs{L'_g} = \delta(L'_g)$. When $F_{g} = \overline{C}_{7}$ (i.e.\ $g = 3$), $c_{g} = 1/3$. We take $L'_g = K_{r - 3}(3) + \overline{C}_{7}$ which is a regular blow-up of $L_{g}$. This satisfies 
	\begin{equation*}
		\delta(L'_g) \abs{L'_g}^{-1} = \frac{3r - 5}{3r - 2} = 1 - \frac{1}{r - 1 + 1/3} = k_{g}.
	\end{equation*}
	This establishes \cref{property:4} except for $g = 5, 7, 9, 11$. \Cref{fig:HHTHweighted} shows the weightings (which induce blow-ups) of $F_{g}$ we would like to use (they maximise the minimum degree relative to the order). However, some care will be needed owing to the zero weights. Consider the case $g = 5$. \Cref{fig:H21weighted} shows a 9-vertex weighted graph with minimum degree 5, however, owing to the zero weights, this is not strictly a blow-up of $F_{5} = H_{2}^{+}$. But for any $k < 5/9$ there is a genuine blow-up, $H_{2}^{+\prime}$, of $H_{2}^{+}$ with $\delta(H_{2}^{+\prime}) \geqslant k \cdot \abs{H_2^{+\prime}}$.
	
	Let $c < k_{5}$: $c = 1 - 1/(r - 1 + \beta)$ for some $\beta < c_{5} = 1/4$ and so $k = 1 - 1/(2 + \beta) < 5/9$. Let $H_{2}^{+\prime}$ be a blow-up of $H_{2}^{+}$ satisfying $\delta(H_{2}^{+\prime}) \geqslant k \cdot \abs{H_2^{+\prime}}$. Let $L_{5}' = K_{r - 3}(\abs{H_2^{+\prime}} - \delta(H_{2}^{+\prime})) + H_{2}^{+\prime}$ which is a blow-up of $L_{5}$. Then
	\begin{align*}
		\delta(L_{5}') \abs{L'_5}^{-1} & = \frac{(r - 3)\abs{H_2^{+\prime}} - (r - 4) \delta(H_{2}^{+\prime})}{(r - 2) \abs{H_2^{+\prime}} - (r - 3) \delta(H_{2}^{+\prime})} \\
		& = 1 - \frac{\abs{H_2^{+\prime}} - \delta(H_{2}^{+\prime})}{(r - 2) \abs{H_2^{+\prime}} - (r - 3) \delta(H_{2}^{+\prime})} \\
		& = 1 - \frac{1}{r - 3 + 1/(1 - \delta(H_{2}^{+\prime}) \abs{H_2^{+\prime}}^{-1})} \\
		& \geqslant 1 - \frac{1}{r - 3 + 1/(1 - k)} \\
		& = 1 - \frac{1}{r - 3 + 2 + \beta} = c,
	\end{align*}
	which establishes \cref{property:4} for $g = 5$. Similar calculations hold for $g = 7, 9, 11$ (for $g = 7$, there are even no troublesome zero weights). Indeed, for $c < k_{g}$: $c = 1 - 1/(r - 1 + \beta)$ for some $\beta < c_{g}$ and so $k = 1 - 1/(2 + \beta)$ is less than 6/11 (for $g = 7$), less than 7/13 (for $g = 9$) and less than 8/15 (for $g = 11$). \Cref{fig:HHTHweighted} shows that there is a blow-up $F_{g}'$ of $F_{g}$ with $\delta(F_{g}') \geqslant k \cdot \abs{F'_g}$. Then $L'_g = K_{r - 3}(\abs{F'_g} - \delta(F_{g}')) + F_{g}'$ is a blow-up of $L_{g}$, which satisfies $\delta(L'_g) \abs{L'_g}^{-1} \geqslant c$. Thus the sequences satisfy \cref{property:4}.
	\begin{figure}[H]
		\centering
		\begin{subfigure}{.24\textwidth}
			\centering
			\begin{tikzpicture}
				\foreach \pt in {0,1,...,6} 
				{
					\tkzDefPoint(\pt*360/7 + 90:1.5){v_\pt}
				} 
				\tkzDefPoint(0,0){u}
				\tkzDrawPolySeg(v_0,v_1,v_2,v_3,v_4,v_5,v_6,v_0)
				\tkzDrawPolySeg(v_1,v_3,v_5,v_0, v_2,v_4,v_6)
				\tkzDrawSegments(u,v_0 u,v_2 u,v_5)
				\tkzDrawPoints(v_0,v_...,v_6)
				\tkzDrawPoint(u)
				\tkzLabelPoint[above](v_0){2}
				\tkzLabelPoint[left](v_1){0}
				\tkzLabelPoint[left](v_2){2}
				\tkzLabelPoint[below](v_3){1}
				\tkzLabelPoint[below](v_4){1}
				\tkzLabelPoint[right](v_5){2}
				\tkzLabelPoint[right](v_6){0}
				\tkzLabelPoint[below](u){1}
			\end{tikzpicture}
			\subcaption{$H_{2}^{+}$ weighted}\label{fig:H21weighted}
		\end{subfigure}
		\begin{subfigure}{.24\textwidth}
			\centering
			\begin{tikzpicture}
				\foreach \pt in {0,1,...,6} 
				{
					\tkzDefPoint(\pt*360/7 + 90:1.5){v_\pt}
				} 
				\tkzDrawPolySeg(v_0,v_1,v_2,v_3,v_4,v_5,v_6,v_0)
				\tkzDrawPolySeg(v_1,v_3,v_5,v_0, v_2,v_4,v_6)
				\tkzDrawPoints(v_0,v_...,v_6)
				\tkzLabelPoint[above](v_0){3}
				\tkzLabelPoint[left](v_1){1}
				\tkzLabelPoint[left](v_2){2}
				\tkzLabelPoint[below](v_3){1}
				\tkzLabelPoint[below](v_4){1}
				\tkzLabelPoint[right](v_5){2}
				\tkzLabelPoint[right](v_6){1}
			\end{tikzpicture}
			\subcaption{$H_{2}$ weighted}\label{fig:H2weighted}
		\end{subfigure}
		\begin{subfigure}{.24\textwidth}
			\centering
			\begin{tikzpicture}
				\foreach \pt in {0,1,...,6} 
				{
					\tkzDefPoint(\pt*360/7 + 90:1.5){v_\pt}
				}
				\tkzDefPoint(0,0.5){t}
				\tkzDefPoint(0.7,-0.6){u_1}
				\tkzDefPoint(-0.7,-0.6){u_6}
				\tkzDrawPolySeg(v_0,v_1,v_2,v_3,v_4,v_5,v_6,v_0) 
				\tkzDrawSegments(t,v_0 t,u_1 t,u_6)
				\foreach \pt in {0,2,3,4,5,6}
				{
					\tkzDrawSegment(u_1,v_\pt)
				}
				\foreach \pt in {0,1,2,3,4,5}
				{
					\tkzDrawSegment(u_6,v_\pt)
				}
				\tkzDrawPoints(v_0,v_...,v_6)
				\tkzDrawPoints(t,u_1,u_6)
				\tkzLabelPoint[above](v_0){4}
				\tkzLabelPoint[left](v_1){0}
				\tkzLabelPoint[left](v_2){0}
				\tkzLabelPoint[below](v_3){1}
				\tkzLabelPoint[below](v_4){1}
				\tkzLabelPoint[right](v_5){0}
				\tkzLabelPoint[right](v_6){0}
				\tkzLabelPoint[below](t){1}
				\tkzLabelPoint[below right = -3pt](u_1){3}
				\tkzLabelPoint[below left = -3pt](u_6){3}
			\end{tikzpicture}
			\subcaption{$T_{0}$ weighted}\label{fig:T0weighted}
		\end{subfigure}
		\begin{subfigure}{.24\textwidth}
			\centering
			\begin{tikzpicture}
				\foreach \pt in {0,1,...,6} 
				{
					\tkzDefPoint(\pt*360/7 + 90:1.5){v_\pt}
				}
				\tkzDefPoint(-0.5,0){ul}
				\tkzDefPoint(0.5,0){ur}
				
				\tkzDrawPolySeg(v_0,v_1,v_2,v_3,v_4,v_5,v_6,v_0) 
				\tkzDrawPolySeg(v_3,v_5,v_0,v_2,v_4)
				\tkzDrawPolySeg(v_6,v_1)
				\tkzDrawSegments(ul,v_0 ul,v_2 ul,v_3 ur,v_0 ur,v_5 ur,v_3)
				\tkzDrawPoints(v_0,v_...,v_6)
				\tkzDrawPoints(ul,ur)
				\tkzLabelPoint[above](v_0){5}
				\tkzLabelPoint[left](v_1){0}
				\tkzLabelPoint[left](v_2){3}
				\tkzLabelPoint[below](v_3){2}
				\tkzLabelPoint[below](v_4){0}
				\tkzLabelPoint[right](v_5){3}
				\tkzLabelPoint[right](v_6){0}
				\tkzLabelPoint[right](ul){1}
				\tkzLabelPoint[left](ur){1}
			\end{tikzpicture}
			\subcaption{$H_{1}^{++}$ weighted}\label{fig:H12weighted}
		\end{subfigure}
		\caption{}\label{fig:HHTHweighted}
	\end{figure}
	
	We are left to show that the sequences satisfy \cref{property:3}. We will make use of the results about $a$-locally bipartite graphs given in \cref{spec4alocalbip}. Let $G$ be a graph with minimum degree greater than $k_{g} \cdot \abs{G}$.
	
	First suppose that $G$ is not $(r - 2)$-locally bipartite. Then there is an $(r - 2)$-clique $K$ in $G$ whose common neighbourhood is not bipartite. Let $G_{K}$ be the induced subgraph of $G$ whose vertex set is the common neighbourhood of the vertices in $K$. Let the vertices of $K$ be $x_{1}, x_{2}, \dotsc, x_{r - 2}$. Note that for each $v \in V(G)$,
	\begin{equation*}
		\mathds{1}(v \in G_{K}) \geqslant \mathds{1}(vx_{1} \in E(G)) + \dotsb + \mathds{1}(vx_{r - 2} \in E(G)) - (r - 3)
	\end{equation*}
	and so summing over all the vertices gives 
	\begin{equation*}
		\abs{G_K} \geqslant (r - 2) \delta(G) - (r - 3) \abs{G} > \bigl((r - 2) k_{g} - (r - 3)\bigr) \cdot \abs{G}.
	\end{equation*}
	Note that $\delta(G_{K}) \geqslant \delta(G) - (\abs{G} - \abs{G_K}) = \abs{G_K} - (\abs{G} - \delta(G))$ so
	\begin{align*}
		\frac{\delta(G_{K})}{\abs{G_K}} & \geqslant 1 - \frac{\abs{G} - \delta(G)}{\abs{G}} \cdot \frac{\abs{G}}{\abs{G_K}} > 1 - \biggl(1 - \frac{\delta(G)}{\abs{G}}\biggr) \cdot \frac{1}{(r - 2) k_{g} - (r - 3)} \\
		& > 1 - \frac{1 - k_{g}}{(r - 2) k_{g} - (r - 3)} = 1 - \frac{1}{1 + c_{g}}.
	\end{align*}
	For $g = 1$: $G_{K}$ is not bipartite and $\delta(G_{K}) > 2/5 \cdot \abs{G_K}$ so, by \cref{deltaHlemma4oddcycles}, $G_{K}$ contains a triangle and so $G$ contains a $K_{r - 2} + K_{3} = L_{0}$. Similarly for $g = 2$, $G$ contains a copy of $L_{0}$ or $L_{1}$ and for $g = 3$, $G$ contains one of $L_{0}$, $L_{1}$ and $L_{2}$. More generally, for $g$ even, $G$ contains $L_{\ell}$ for some $\ell \in \set{0, 1, 2, 4, 6, \dotsc, g - 2}$ and for $g$ odd, $G$ contains $L_{\ell}$ for some $\ell \in \set{0, 1, 2, 4, 6, \dotsc, g - 1}$.
	
	Now suppose that $G$ is $(r - 2)$-locally bipartite. For $g = 1, 2, 3$:
	\begin{equation*}
		\delta(G)\abs{G}^{-1} > 1 - \frac{1}{r - 1 + c_{3}} = 1 - \frac{1}{r - 2 + 4/3}, 
	\end{equation*}
	so, by \cref{spec4alocalbip}, $G$ is $r$-colourable. For $g = 4, 5$:
	\begin{equation*}
		\delta(G)\abs{G}^{-1} > 1 - \frac{1}{r - 1 + c_{5}} = 1 - \frac{1}{r - 2 + 5/4},
	\end{equation*}
	so, by \cref{spec4alocalbip}, $G$ is $r$-colourable or contains $K_{r - 3} + \overline{C}_{7} = L_{3}$. For $g = 6, 7$:
	\begin{equation*}
		\delta(G)\abs{G}^{-1} > 1 - \frac{1}{r - 1 + c_{7}} = 1 - \frac{1}{r - 2 + 6/5},
	\end{equation*}
	so, by \cref{spec4alocalbip}, $G$ is $r$-colourable or contains $K_{r - 3} + \overline{C}_{7} = L_{3}$ or contains $K_{r - 3} + H_{2}^{+} = L_{5}$. For $g = 8, 9$:
	\begin{equation*}
		\delta(G)\abs{G}^{-1} > 1 - \frac{1}{r - 1 + c_{9}} = 1 - \frac{1}{r - 2 + 7/6},
	\end{equation*}
	so, by \cref{spec4alocalbip}, $G$ is $r$-colourable or contains $K_{r - 3} + H_{2} = L_{7}$. Finally, for $g = 10, 11$:
	\begin{equation*}
		\delta(G)\abs{G}^{-1} > 1 - \frac{1}{r - 1 + c_{11}} = 1 - \frac{1}{r - 2 + 8/7},
	\end{equation*}
	so, by \cref{spec4alocalbip}, $G$ is $r$-colourable or contains $K_{r - 3} + H_{2} = L_{7}$ or contains $K_{r - 3} + T_{0} = L_{9}$. Hence, the sequences do indeed satisfy \cref{property:3}, as required.
\end{Proof}

\section{Approximate chromatic profile and homomorphism-free graphs}\label{sec:approx}

In \cref{sec:context} we introduced two chromatic profiles related to $\delta_{H}$. The first was the approximate chromatic profile which for a family $\cF$ is,
\begin{equation*}
	\begin{split}
		\delta^{\ast}_{\chi}(\cF, k) = \inf \set{c \colon &\textnormal{if }\abs{G} = n, \delta(G) \geqslant c n, G \in \cF, \\
		& \textnormal{then } G \textnormal{ can be made } k\textnormal{-colourable by deleting at most } o(n^{2}) \textnormal{ edges}},
	\end{split}
\end{equation*}
and the second was the chromatic profile of the family of $H$-homomorphism-free graphs,
\begin{equation*}
	\delta_{\chi}(H\text{-hom}, k) = \inf\set{c \colon \text{if }\delta(G) \geqslant c \abs{G} \text{ and } H \nrightarrow G, \text{then } G \text{ is $k$-colourable}},
\end{equation*}
We first prove \cref{profiles}, showing that families which are closed under blow-ups have nicely behaved chromatic profiles.

\profiles*

\begin{Proof}
	Plainly the definitions of the approximate and normal chromatic profiles give
	\begin{equation*}
		\delta_{\chi}(\cF, k) \geqslant \delta_{\chi}^{\ast}(\cF, k),
	\end{equation*}
	as being $k$-colourable implies being within $o(n^{2})$ edges of $k$-colourable. 
	
	Let $c < \delta_{\chi}(\cF, k)$. By definition, there must be a graph $G \in \cF$ which is not $k$-colourable and has minimum degree at least $c \abs{G}$. Let $G'$ be the balanced blow-up of $G$ on $n$ vertices: as $\cF$ is closed under taking blow-ups, $G' \in \cF$. Furthermore, $G'$ has minimum degree at least $c \abs{G} \floor{n/\abs{G}} = (c - o(1)) n$. Finally, just as in \cref{sec:deltaHprop}, it is easy to check by counting copies of $G$ that to make $G'$ $k$-colourable requires the deletion of $\Omega(n^{2})$ edges. In particular, \(\delta_{\chi}^{\ast}(\cF, k) \geqslant c\), as required.
\end{Proof}

Finally, we prove \cref{profileshom} showing that the chromatic profile of $H$-homomorphism-free graphs and the approximate chromatic profile of $H$-free graphs are one and the same.

\profileshom*

\begin{Proof}
	We first show that $\delta^{\ast}_{\chi}(H, k) \geqslant \delta_{\chi}(H\text{-hom}, k)$. Any graph to which there is no homomorphism from $H$ must be $H$-free, so $\delta^{\ast}_{\chi}(H, k) \geqslant \delta^{\ast}_{\chi}(H\text{-hom}, k)$. But $H$-hom is closed under taking blow-ups so, by \cref{profiles}, $\delta^{\ast}_{\chi}(H\text{-hom}, k) = \delta_{\chi}(H\text{-hom}, k)$, as required.
	
	We now show that $\delta^{\ast}_{\chi}(H, k) \leqslant \delta_{\chi}(H\text{-hom}, k)$. Let $\gamma > 0$, $c = \delta_{\chi}(H\text{-hom}, k) + 2\gamma$ and let $G$ be an $n$-vertex $H$-free graph with $\delta(G) \geqslant c n$. Let 
	\begin{equation*}
		\cH = \set{H' \colon H \to H', \abs{H'} \leqslant \abs{H}},
	\end{equation*}
	and note that this is finite and that there is a homomorphism from $H$ to a graph if and only if that graph contains some $H' \in \cH$. There is some $t \leqslant \abs{H}$ such that $H \subset H'(t)$ for every $H' \in \cH$. 
	
	Fix $H' \in \cH$. The graph $G$ is $H$-free so does not contain $H'(t)$. By Erd\H{o}s's result on the extremal function for complete $\ell$-uniform $\ell$-partite hypergraphs~\cite{Erdos1964hypextremal}, $G$ must contain $o(n^{\abs{H'}})$ copies of $H'$ (see for example~\cite[Lemma 6.2]{ABGKM2017dense}). By the graph removal lemma (see for example~\cite[Theorem 2.9]{KomlosSimonovits1996}), $G$ can be made $H'$-free by deleting $o(n^{2})$ edges. As $\cH$ is finite, there is a spanning subgraph $G'$ of $G$ with $e(G) - e(G') = o(n^{2})$ which contains no $H' \in \cH$.
	
	Take $G'$ and sequentially delete vertices of degree less than $(c - \gamma)n$ until no more remain. Provided $n$ is large enough, so that $(e(G) - e(G'))/n^{2}$ is sufficiently small, this process will terminate with the deletion of at most $o(n)$ vertices. Let the resulting graph be $G''$. Then $G''$ satisfies $\delta(G'') \geqslant (c - \gamma) \abs{G''}$ and $G''$ contains no $H' \in \cH$. In particular, there is no homomorphism $H \to G''$. But $c - \gamma > \delta_{\chi}(H\text{-hom}, k)$ so $G''$ is $k$-colourable. Furthermore $G''$ was obtained from $G$ by the deletion of $o(n^{2})$ edges. Thus $\delta^{\ast}_{\chi}(H, k) \leqslant c$, as required.
\end{Proof}

\section{Concluding remarks}

While \cref{deltaH4} determines $\delta_{H}$ for most $H$, many minimum degree stability questions remain. Firstly how do those $\delta_{H}$ not determined by the theorem behave? Is $1 - 1/(r - 1)$ the only accumulation point of $\set{\delta_{H} \colon \chi(H) = r + 1}$ (as it is for $r = 2$) or is there more exotic behaviour? To extend \cref{deltaH4} one would need to extend \cref{spec4alocalbip} below $1 - 1/(a + 8/7)$ and so we would need to further our knowledge of the structure of locally bipartite and more generally locally colourable graphs -- for more information see \cite{Illingworth2022localbpart}. All the graphs appearing in \cref{deltaH4} are either an $(r - 2)$-clique joined to an odd cycle or an $(r - 3)$-clique joined to a 4-chromatic locally bipartite graphs. Are these the only graphs that appear? The motivation in \cref{sec:graphs} suggests that this is the case for $r = 3$: the major obstacles to being close to tripartite are containing either a blow-up of an odd wheel of the blow-up of some 4-chromatic locally bipartite graph. However, for greater $r$ other graphs could appear. For example, containing a blow-up of some 5-chromatic locally tripartite graph could be an obstacle for being close to 4-partite.

In \cref{sec:context} we placed $\delta_{H}$ as the first non-trivial threshold within the approximate chromatic profile. It would be interesting to determine the next threshold, that is, to understand the behaviour of $\delta_{\chi}^{\ast}(H, \chi(H))$. For triangles and cliques, this has already been done~\cite{BrandtThomasse2005,GoddardLyle2010,Nikiforov2010}.

Another natural direction is to consider the number of edges that need deleting to make an $n$-vertex $H$-free graph, $G$, with minimum degree at least $(\delta_{H} + \varepsilon) n$ $r$-partite. Are $o(n^{2})$ edges really required, or, as is often the case, can one get away with $\cO(n^{2 - \rho})$ for some $\rho = \rho(H) > 0$? This has precedent. Erd\H{o}s and Simonovits~\cite{Erdos1967, Erdos1968, Simonovits1968} showed that the $H$-free graph with most edges can be made $r$-partite by deleting $\cO(n^{2 - \rho})$ edges. Also Alon and Sudakov's result, \cref{AS}, gives an affirmative answer for $H = K_{r + 1}(t)$. The heuristic here is that if more than $n^{2 - \rho}$ edges are required, then, by the theorem of K\H{o}v\'{a}ri, S\'{o}s and Tur\'{a}n~\cite{KST1954}, $G$ is an $r$-partite graph with some large $K_{t, t}$ appearing inside one of the parts. Joining these together ought to give some blow-up of an $F_{g}$ and so a copy of $H$. In both the cases of Erd\H{o}s-Simonovits and Alon-Sudakov, the minimum degree of $G$ was large and so $G$ was well connected. For our present situation the following would be the most basic question. Although an affirmative answer seems plausible (and I have a proof for $g = 3$) the smaller minimum degree casts doubt for large $g$.
\begin{question}
	For positive integers $g$ and $t$ is there some $\rho > 0$ such that every $n$-vertex graph with minimum degree at least $(2/(2g + 1) + \varepsilon) n$ either contains $C_{2g - 1}(t)$ or can be made bipartite by deleting $\cO(n^{2 - \rho})$ edges?
\end{question}

\section*{Acknowledgements}

It is a pleasure to thank Andrew Thomason for many helpful discussions. I am grateful to the anonymous referees for their careful reading and excellent suggestions for improving the presentation.

{
\fontsize{11pt}{12pt}
\selectfont
	
\hypersetup{linkcolor={red!70!black}}
\setlength{\parskip}{2pt plus 0.3ex minus 0.3ex}
	
\newcommand{\etalchar}[1]{$^{#1}$}

}

\end{document}